\documentclass[a4paper]{amsart}

\usepackage{amsmath,amssymb,amsthm}
\usepackage{fullpage}
\usepackage{natbib}
\usepackage{graphicx}
\usepackage{caption}

\newcommand{\R}{\ensuremath{\mathbb{R}}}

\begin{document}

\title{Topology Applied to Machine Learning:\\ From Global to Local} 
\author{Henry Adams}
\email{henry.adams@colostate.edu}
\author{Michael Moy}
\email{michael.moy@colostate.edu}

\maketitle

\begin{abstract}
Through the use of examples, we explain one way in which applied topology has evolved since the birth of persistent homology in the early 2000s.
The first applications of topology to data emphasized the \emph{global} shape of a dataset, such as the three-circle model for $3 \times 3$ pixel patches from natural images, or the configuration space of the cyclo-octane molecule, which is a sphere with a Klein bottle attached via two circles of singularity.
In these studies of global shape, short persistent homology bars are disregarded as sampling noise.
More recently, however, persistent homology has been used to address questions about the \emph{local} geometry of data.
For instance, how can local geometry be vectorized for use in machine learning problems?
Persistent homology and its vectorization methods, including persistence landscapes and persistence images, provide popular techniques for incorporating both local geometry and global topology into machine learning.
Our meta-hypothesis is that the short bars are as important as the long bars for many machine learning tasks.
In defense of this claim, we survey applications of persistent homology to shape recognition, agent-based modeling, materials science, archaeology, and biology.  
Additionally, we survey work connecting persistent homology to geometric features of spaces, including curvature and fractal dimension, and various methods that have been used to incorporate persistent homology into machine learning.


\end{abstract}

\section{Introduction}

Applied topology is designed to measure the \emph{shape} of data --- but what is shape?
Early examples in applied topology found low-dimensional structures in high-dimensional datasets, such as the three circle and Klein bottle models for greyscale natural image patches.
These models are global: they parameterize the entire dataset, in the sense that most of the data points look like some point in the model, plus noise.
In more recent applications, however, the shape that is being measured is not global, but instead local.
Local features include texture, small-scale geometry, and the structure of noise.

Indeed, for the first decade after the invention of persistent homology, the primary story was that significant features in a dataset corresponded to long bars in the persistence barcode, whereas shorter bars generally corresponded to sampling noise.
This story has evolved as applied topology has become incorporated into the machine learning pipeline.
In machine learning applications, many researchers have independently found (as we survey in Sections~\ref{sec:local}--\ref{sec:ML}) that the short bars are often the most discriminating --- the shape of the noise, or of the local geometry, is what often enables high classification accuracy.
We want to emphasize that short bars do matter.
Indeed, the short bars in persistent homology are currently one of the best out-of-the-box methods for summarizing local geometry for use in machine learning.
Though humans may not be able to interpret short persistent homology bars on our own (there may be too many short bars for the human eye to count), machine learning algorithms can be trained to do so.
In this way, persistent homology has greatly expanded in scope during the second decade after its invention: persistent homology has important applications as a descriptor not only of global shape, but also of local geometry.

In this perspective article, we begin by outlining some of the most famous early applications of persistent homology in the global analysis of data, in which short bars were disregarded as noise.
Our meta-hypothesis, however, is that short bars do matter, and furthermore, they matter crucially when combining topology with machine learning.
As a partial defense for this claim, we provide a selected survey on the use of persistent homology in measuring texture, noise, local geometry, fractal dimension, and local curvature.
We predict that the applications of persistent homology to machine learning will continue to advance in number, impact, and scope, as persistent homology is a mathematically motivated out-of-the-box tool that one can use to summarize not only the global topology but also the local geometry of a wide variety of datasets.

\section{Point cloud and sublevel set persistent homology}

What is a persistent homology bar?
The homology of a space, roughly speaking, records how many holes that a space has in each dimension.
A 0-dimensional hole is a connected component, a 1-dimensional hole is a loop, a 2-dimensional hole is a void enclosed by a surface like a sphere or a torus, etc.
Homology becomes persistent when one is instead given a \emph{filtration}, i.e.\ an increasing sequence of spaces.
Each hole is now represented by a bar, where the start (resp.\ end) point of the bar corresponds to the first (resp.\ last) stage in the filtration where the topological feature is present~\citep{edelsbrunner2000topological}.
Short bars correspond to features with short lifetimes, which are quickly filled-in after being created.
By contrast, long bars correspond to more \emph{persistent} features.

\begin{figure}[h]
\captionsetup{width=0.95\textwidth}
\begin{center}
\includegraphics[width=\textwidth]{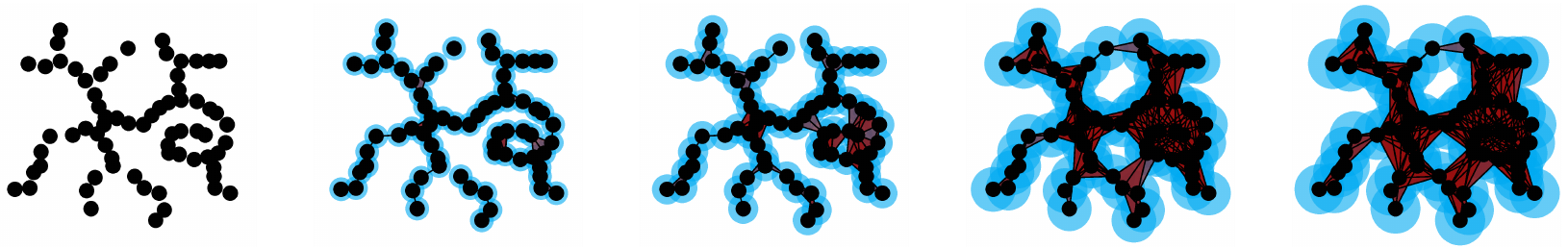}
\caption{A point cloud, the surrounding union of balls, and its \v{C}ech complexes at different choices of scale.}
\label{fig:CechExample}
\end{center}
\end{figure}

Perhaps the two most frequent contexts in which persistent homology is applied are point cloud persistent homology and sublevel set persistent homology.
In \emph{point cloud persistent homology}, the input is a finite set of points (a point cloud) residing in Euclidean space or some other metric space~\citep{Carlsson2009}.
For any real number $r>0$, we consider the union of all balls of radius $r$ centered at some point in our point cloud; see Figure~\ref{fig:CechExample}.
This union of balls provides our filtration as the radius $r$ increases.\footnote{In practice, the union of balls is stored or approximated by a \emph{simplicial complex}, for example a \emph{\v{C}ech} or \emph{Vietoris--Rips} complex~\citep{ChazalDeSilvaOudot2013}.}
A typical interpretation of the resulting persistent homology, from the global perspective, is that the long persistent homology bars recover the homology of the ``true" underlying space from which the point cloud was sampled~\citep{ChazalOudot2008}.
A more modern but increasingly utilized perspective is that the short persistent homology bars recover the local geometry --- i.e.\ the texture, curvature, or fractal dimension of the point cloud data.

\begin{figure}[h]
\captionsetup{width=0.95\textwidth}
\begin{center}
\includegraphics[width=\textwidth]{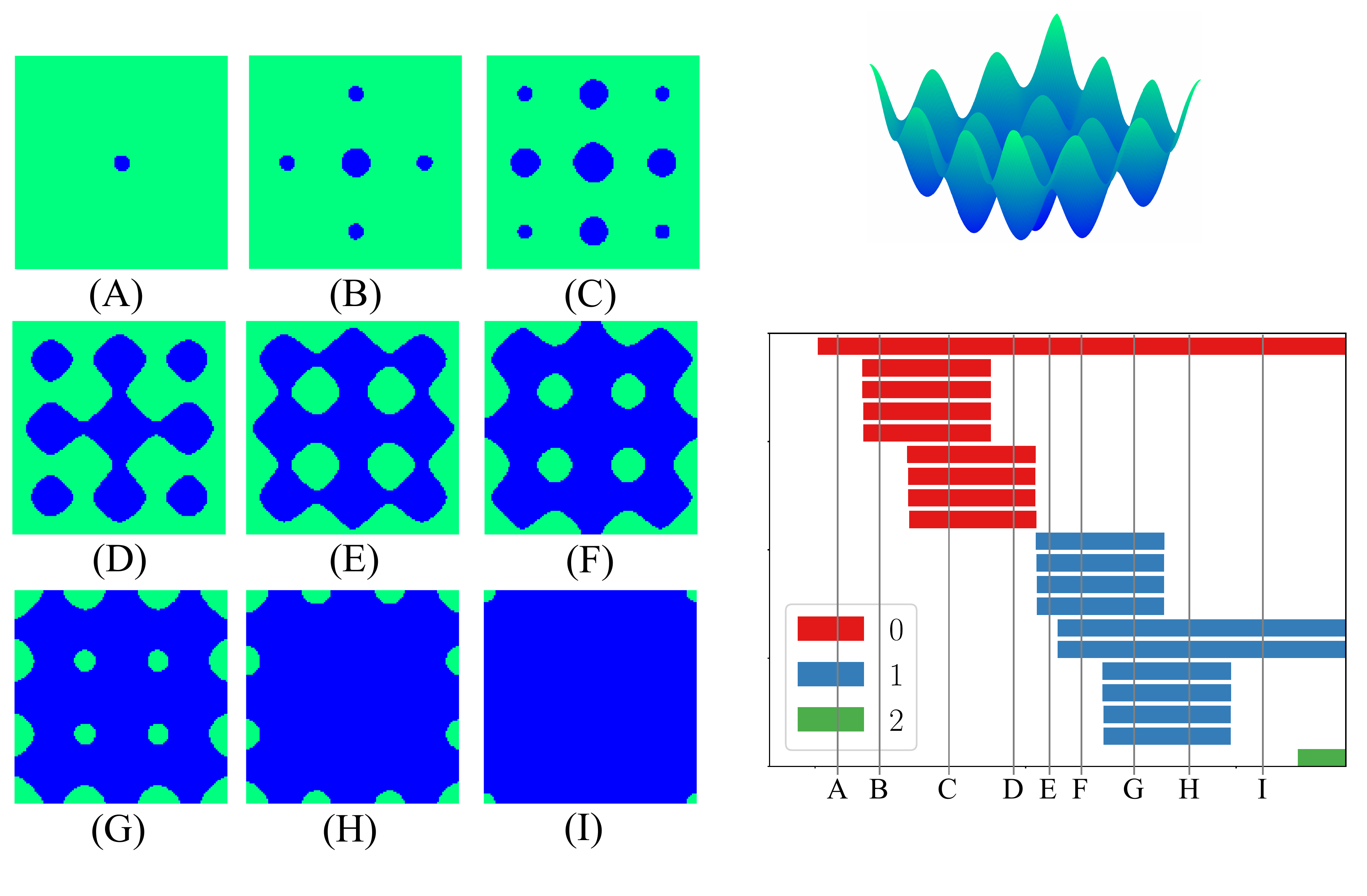}
\caption{(Top right)
An energy function for the molecule pentane.
The domain is a torus, i.e.\ a square with periodic boundary conditions, as there are two circular degrees of freedom (dihedral angles) in the molecule. (Left) Nine different sublevel sets of energy. (Bottom right) The sublevel set persistent homology of this energy function on the torus, with 0-dimensional homology in red, 1-dimensional homology in blue, 2-dimensional homlogy in green.
Image from~\cite{mirth2020representations}.}
\label{fig:sublevelset-ph-pentane}
\end{center}
\end{figure}

In \emph{sublevel set persistent homology}, the input is instead a real-valued function $f\colon Y\to\R$ defined on a space $Y$~\citep{cohen2007stability}.
For example, $Y$ may be a Euclidean space of some dimension.
The filtration arises by considering the sublevel sets $\{y\in Y~|~f(y)\le r\}$.
As the threshold $r$ increases, the sublevel sets grow.
One can think of $f$ as encoding an \emph{energy}, in which case sublevel set persistent homology encodes the shape of low-energy configurations~\citep{mirth2020representations}.
The length of a bar then measures how large of an energy barrier must be exceeded in order for a topological feature to be filled-in: a short bar corresponds to a feature that is quickly filled-in by exceeding a low energy barrier, whereas a long bar corresponds to a topological feature that persists over a longer range of energies; see Figure~\ref{fig:sublevelset-ph-pentane}.
Sublevel set persistent homology is frequently applied to grayscale image data or matrix data, where a real-valued entry of the image or matrix is interpreted as the value of the function $f$ on a pixel.

We remark that the ``union of balls" filtration for point cloud persistent homology can be viewed as a version of sublevel set persistent homology: a union of balls of radius $r$ is the sublevel set at threshold $r$ of the distance function to the set of points in the point cloud.

\begin{figure}[h!]
\begin{center}
\includegraphics[width=5in]{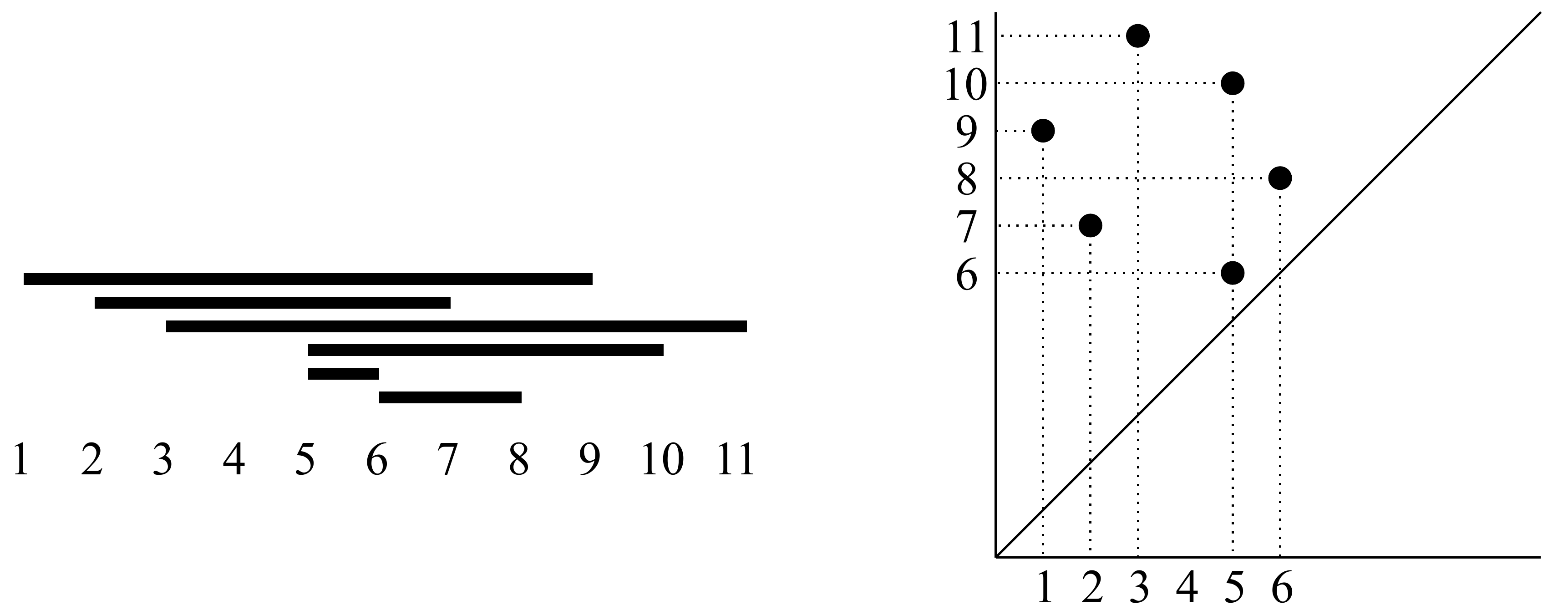}
\end{center}
\caption{(Left) A persistent homology barcode, with the birth and death scale of each bar indicated on the horizontal axis.
(Right) Its corresponding persistence diagram, i.e.\ a collection of points in the first quadrant above the diagonal, with birth coordinates on the horizontal axis and death coordinates on the vertical axis.}
\label{fig:barcode-diagram}
\end{figure}

Persistent homology can be represented in two equivalent ways: either as a persistence barcode or as a persistence diagram; see Figure~\ref{fig:barcode-diagram}.
Each interval in the persistence barcode is represented in the persistence diagram by a point in the plane, with its birth coordinate on the horizontal axis and with its death coordinate on the vertical axis\footnote{Barcodes allow for open or closed endpoints of intervals.
This information can be also be encoded in a \emph{decorated} persistence diagram~\citep{chazal2016structure}.}.
As the death of each feature is after its birth, persistence diagram points all lie above the diagonal line $y=x$.
Short bars in the barcode correspond to persistence diagram points close to the diagonal, and long bars in the barcode correspond to persistence diagram points far from the diagonal.

\section{Examples measuring global shape}

The earliest applications of topology to data measured the global shape of a dataset.
In these examples, the long persistent homology bars represented the true homology underlying the data, whereas the small bars were ignored as artifacts of sampling noise.

\begin{figure}[h!]
\begin{center}
\includegraphics[width=3in]{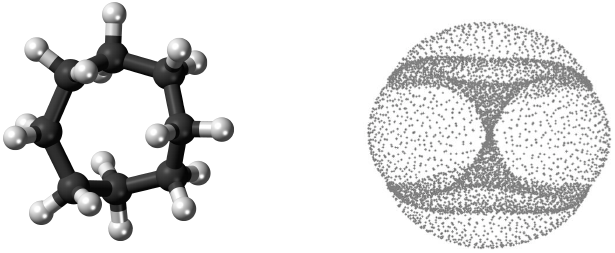}
\end{center}
\caption{(Left) The cyclo-octane molecule consists of a ring of 8 carbon atoms (black), each bonded to a pair of hydrogen atoms (white).
(Right) A PCA projection of a dataset of different conformations of the cyclo-octane molecule; this shape is a sphere glued to a Klein bottle (the ``hourglass") along two circles of singularity.
The right image is from~\citet{martin2010topology}.}
\label{fig:cyclo}
\end{figure}

What do we mean by ``global shape"?
Consider, for example, conformations of the cyclo-octane molecule $C_8H_{16}$, which consists of a ring of eight carbons atoms, each bonded to a pair of hydrogen atoms; see Figure~\ref{fig:cyclo} (left).
The locations of the carbon atoms in a conformation approximately determine the locations of the hydrogen atoms via energy minimization, and hence each molecule conformation can be mapped to a point in $\R^{24}=\R^{8\cdot 3}$, as the location of each carbon atom can be specified by three coordinates.
This map realizes the conformation space of cyclo-octane as a subset of $\R^{24}$, and then we mod out by rigid rotations and translations.
Topologically, the conformation space of cyclo-octane turns out to be the union of a sphere with a Klein bottle, glued together along two circles of singularities; see Figure~\ref{fig:cyclo} (right).  
This model was obtained by~\citet{martin2010topology,martin2011non,brown2008algorithmic}, who furthermore obtain a triangulation of this dataset (a representation of the dataset as a union of vertices, edges, and triangles).

\begin{figure}[h!]
\begin{center}
\includegraphics[width=6in]{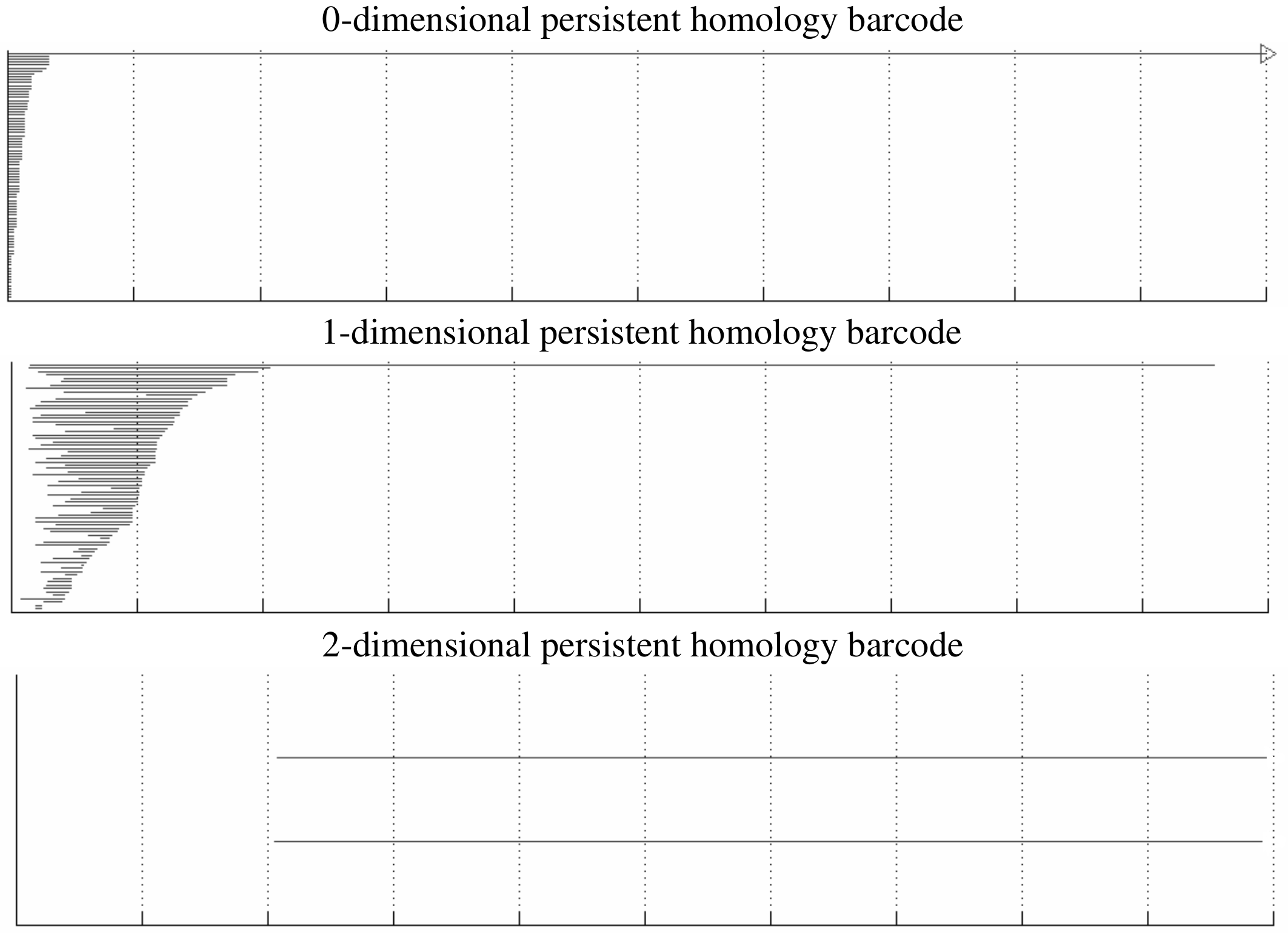}
\end{center}
\caption{0-, 1-, and 2-dimensional persistent homology barcodes for the cyclo-octane dataset.
The horizontal axis corresponds to the birth and death scale of the bars, and the vertical axis is an arbitrary ordering of the bars (here by death scale).}
\label{fig:cyclo-barcodes}
\end{figure}

A Klein bottle, like a sphere, is a 2-dimensional manifold.
Whereas a sphere can be embedded in 3-dimensional space, a Klein bottle requires at least four dimensions in order to be embedded without self-intersections.
When a sphere and Klein bottle are glued together along two circles, the union is no longer a manifold.
Indeed, near the gluing circles, the space does not look like a sheet of paper, but instead like the tail of a dart with four fins, i.e., the letter `X' crossed with the interval $[0,1]$.
However, the result is still a 2-dimensional stratified space.
In Figure~\ref{fig:cyclo-barcodes}, we compute the persistent homology of a point cloud dataset of 1,000,000 cyclo-octane molecule configurations.
The short bars are interpreted as noise, whereas the long bars are interpreted as attributes of the underlying shape.
We obtain a single connected component, a single 1-dimensional hole, and two 2-dimensional homology features.
These homology signatures agree with the homology of the union of a sphere with a Klein bottle, glued together along two circles of singularities.

One of the first applications of persistent homology was to measure the global shape of a dataset of image patches~\citep{CarlssonIshkhanovDeSilvaZomorodian2008}.
This dataset of natural $3\times 3$ pixel patches from black-and-white photographs from indoor and outdoor scenes in fact has three different global shapes!
The most common patches lie along a circle of possible directions of linear gradient patches (varying from black to gray to white).
The next most common patches lie along a three circle model, additionally including a circle's worth of horizontal quadratic gradients, and a circle's worth of vertical quadratic gradients.
At the next level of resolution, the most common patches in some sense lie along a Klein bottle.
All three of these models --- the circle, the three circles, and the Klein bottle --- are global models, summarizing the global shape of the dataset at different resolutions.

\section{Examples measuring local geometry}
\label{sec:local}

Though a single long bar in persistent homology may carry a lot of information, a single small bar typically does not.
However, together a collection of small bars may unexpectedly carry a large amount of geometric content.
A long bar is a trumpet solo --- piercing through to be heard over the orchestra with ease.
The small bars are the string section --- each small bar on its own is relatively quiet, but in concert the small bars together deliver a powerful message.
We survey several modern examples where small persistent homology bars are now the signal, instead of the noise.

\begin{figure}[h!]
\begin{center}
\includegraphics[width=3.5in]{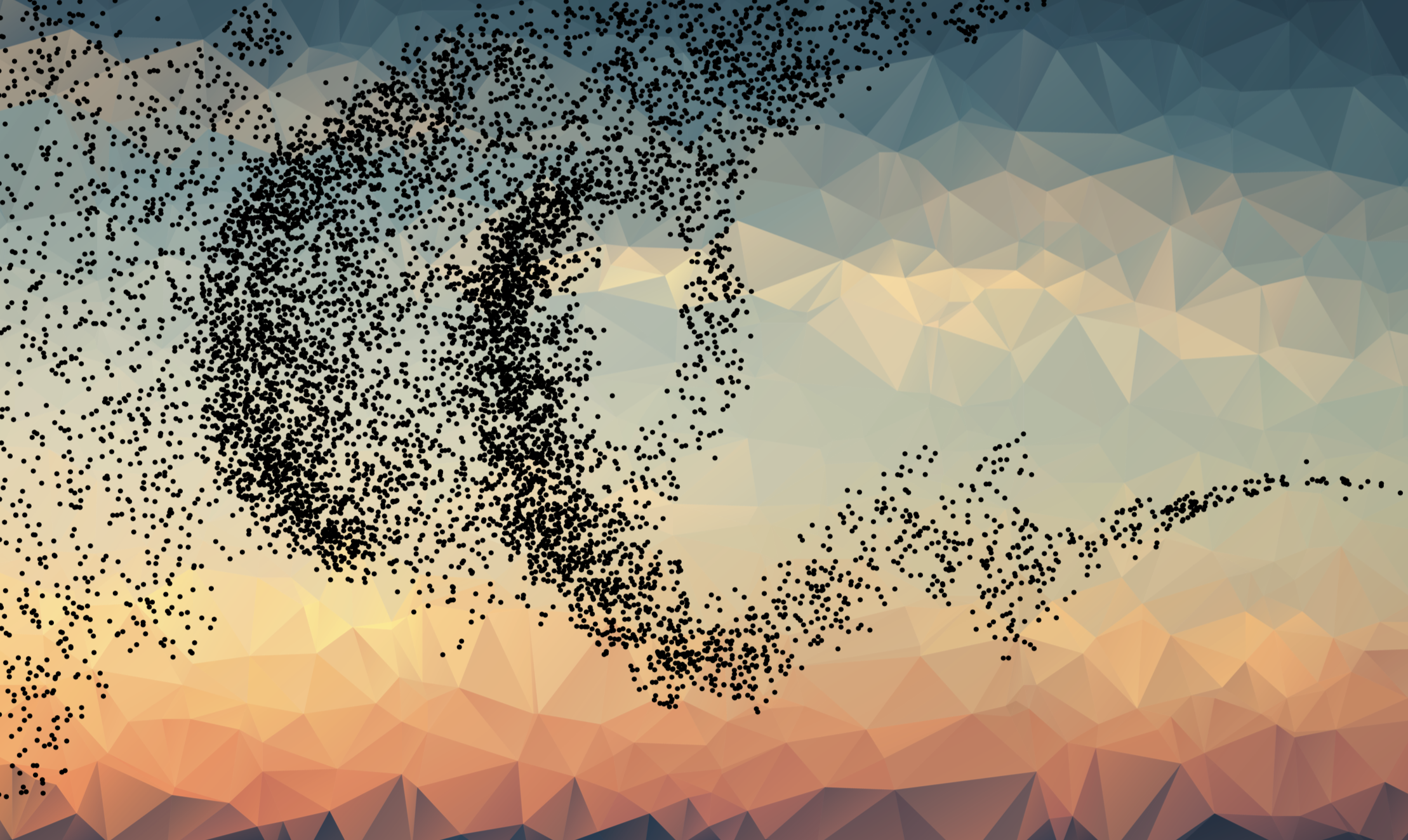}
\end{center}
\caption{A large amount of local and global geometric information is contained in a flock of birds.}
\label{fig:flock-of-birds}
\end{figure}

Birds, fish, and insects move as flocks, schools, and hordes in a way which is determined by \emph{collective motion}: each animal's next motion is a random function of the location of its nearby neighbors.
In a flock of thousands of birds, there is an impressively large amount of time-varying geometry, including for example all $\binom{n}{2}$ pairwise distances, where $n$ is the number of birds; see Figure~\ref{fig:flock-of-birds}.
How can one summarize this much geometric content for use in machine learning tasks, say to predict how the motion of the flock will vary next, or to predict some of the parameters in a mathematical model approximately governing the motion of the birds?
Persistent homology has been used in~\citet{topaz2015topological,ulmer2019topological,bhaskar2019analyzing,TDAcollective,xian2020capturing} to reduce a large collection of geometric content down to a concise summary.
These datasets of animal swarms do not lie along beautiful manifolds (global shapes), but nevertheless there is a wealth of information in the local geometry as measured by the short persistent homology bars.
For example,~\cite{ulmer2019topological} show via time-varying persistent homology\footnote{In particular, the \emph{crocker plot}~\citep{topaz2015topological}} that a control model for aphid motion, in which aphids move independently at random, does not fit experimental data as well as a model incorporating social interaction (distances to nearby neighbors) between the aphids.

Other recent work has used persistent homology to characterize the complexity of geometric objects.
\citet{bendich2016persistent} apply sublevel set persistent homology to the study of brain artery trees, examining the effects of age and sex on the barcodes generated from artery trees.
While younger brains have artery trees containing more local twisting and branching, older brains are sparser with fewer small branches and leaves.
The authors use the 100 longest bars in dimensions 0 and 1 in their analysis, and they further examine which lengths of bars give the highest correlation with age and sex.
For instance, when examining age, they find it is not the longest bars, but instead the bars of medium length (roughly the 21st through 40th longest bars) that are the most discriminatory.

\begin{figure}[h!]
\begin{center}
\includegraphics[width=5in]{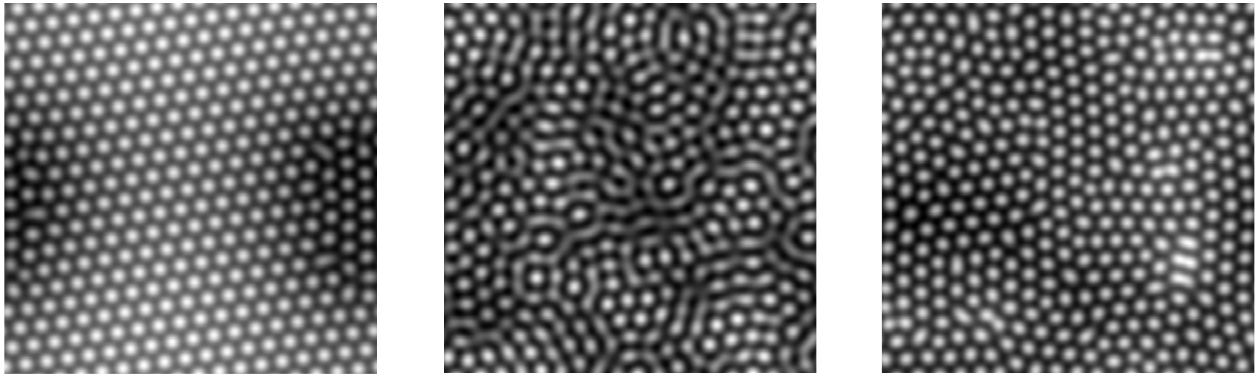}
\end{center}
\caption{Hexagonal lattices, of varying degree of regularity, created by ion bombardment.
Figures from~\cite{motta2018measures}.}
\label{fig:hexagonal}
\end{figure}

In other datasets where points are nearly evenly spaced, barcodes will consist of bars with mostly similar birth and death times.  
Consider for instance the point cloud persistent homology for a square grid of points in the plane: all 0-dimensional bars are identical and adding a small amount of noise to the points will result in a small change to the bars.
The same is true for 1-dimensional bars.  
With this in mind,~\citet{motta2018measures} use persistent homology to measure the order, or regularity, of lattice-like datasets, focusing on hexagonal grids formed by ion bombardment of solid surfaces; see Figure~\ref{fig:hexagonal}. 
The authors' techniques use the variance of 0-dimensional homology bar lengths, and the sum of the lengths of 1-dimensional homology bars, as well as a particular linear combination of the two especially suited to hexagonal lattices.  
Their results suggest that techniques based on persistent homology can provide useful measures of order that are sensitive to both large scale and small scale defects in lattices.
Point cloud persistence has also been used to summarize the local order and randomness in other materials science and chemistry contexts, including amorphous solids and glass~\citep{hiraoka2016hierarchical,nakamura2015persistent,hirata2020structural}, nanoporous materials used in gas adsorption~\citep{krishnapriyan2020persistent}, crystal structure~\citep{Maroulas2020Bayesian}, and protein folding~\citep{xia2014persistent,cang2018integration}.

Though the above examples focus on point cloud persistence, sublevel set persistent homology has also been used to detect the local geometry of functions.
\citet{kramar2016analysis} use sublevel set persistence to summarize the complicated spatio-temporal patterns that arise from dynamical systems modeling fluid flow, including turbulence (Kolmogorov flow) and heat convection (Rayleigh-B\'{e}nard convection).
With sublevel set persistence,~\citet{zeppelzauer2016topological} improve 3D surface classification, including on an archaeology task of segmenting engraved regions of rock from the surrounding natural rock surface.
In a task of tracking automobiles, \citet{bendich2016topological} use the sublevel set persistent homology of driver speeds in order to characterize driver behaviors and prune out improbable paths from their multiple hypothesis tracking framework.



\section{Theory of how persistent homology measures local geometry}

Recent work has begun to formalize the idea that persistent homology measures local geometry.  
\citet{bubenik2020persistent} explore the effect of the curvature of a space on the persistent homology of a sample of points, focusing on disks in spaces with constant curvature.  
Their work includes theoretical results about the persistence of triangles in these spaces, and they are also able to demonstrate experimentally that persistent homology in dimensions 0 and 1 can be used to accurately estimate the curvature given a random sample of points.  
Since the disks in spaces with different curvature are homeomorphic, the differences in persistent homology cannot be due to topology, but rather result from the geometric features of the spaces.

Fractal dimension is another measure of local geometry, and indeed some of the earliest applications of persistent homology in Vanessa Robins' PhD thesis were motivated as a way to capture the fractal dimension of an infinite set in Euclidean space~\citep{robins2000computational,macpherson2012measuring}.
Can this also be applied to datasets, i.e.\ to random collections of finite sets of points?
Given a random sample of points from a measure,~\citet{adams2018fractal} use persistent homology to detect the fractal dimension of the support of the measure.
This notion of \emph{persistent homology fractal dimension} agrees with the Hausdorff/box-counting dimension for 0-dimensional persistent homology and a restricted class of measures; see~\citet{schweinhart2019persistent,schweinhart2020fractal} for further theoretical developments.

A related line of work studies what can be proven about the topology of random point clouds, typically as the number of points in the point cloud goes to infinity~\citep{adler2014crackle,bobrowski2014topology,bobrowski2015maximally,kahle2011random}.
The magnitude~\citep{leinster2013magnitude} and magnitude homology~\citep{hepworth2017categorifying,leinster2017magnitude} of a metric space measure both local and global properties; recent and ongoing work is being done to connect magnitude with persistent homology~\citep{otter2018magnitude,govc2021persistent}.
See also~\citet{weinberger2019interpolation} for connections between sublevel set persistent homology and the geometry of spaces of functions, including Lipschitz constants of functions.
We predict that much more work demonstrating how local geometric features can be recovered from persistent homology barcodes will take place over the next decade.

\section{Machine learning}
\label{sec:ML}

Because persistent homology gives a concise description of the shape of data, it is not surprising that recent work has incorporated persistent homology into machine learning.
When might one consider using persistent homology in concert with machine learning, as opposed to other more classical machine learning techniques measuring shape such as clustering~\citep{xu2005survey} or nonlinear dimensionality reduction~\citep{tenenbaum2000global,roweis2000nonlinear,mcinnes2018umap,kohonen2012self}?
We recommend persistent homology when one desires either (i) a quantitative reductive summary of local geometry, (ii) an estimate of the number or size of more global topological features in a dataset, or (iii) a way to explore if either local geometry or global topology may be discriminatory for the machine learning task at hand.
Researchers have taken at least three distinct approaches: persistence barcodes have been adapted to be input to machine learning algorithms, topological methods have been used to create new algorithms, and persistent homology has been used to analyze machine learning algorithms.

Perhaps the most natural of these approaches is inputting persistence data into a machine learning algorithm.
Though the persistent homology bars provide a summary of both local geometry and global topology, for a quantitative summary to be fully applicable it needs to be amenable for use in machine learning tasks.
The space of persistence barcodes is not immediately appropriate for machine learning.
Indeed, averages of barcodes need not be unique~\citep{mileyko2011probability}, and the space of persistence barcodes does not coarsely embed into any Hilbert space~\citep{bubenik2020embeddings}.
These limitations have initiated a large amount of research on transforming persistence barcodes into more natural formats for machine learning.
From barcodes,~\citet{bubenik2015statistical} creates \emph{persistence landscapes}, which live in a Banach space of functions\footnote{In practice, a different metric is sometimes chosen to map landscapes into a Hilbert space, though the restrictions of~\citet{bubenik2020embeddings} apply.}.
Persistence landscapes are created by rotating a persistence diagram on its side --- so that the diagonal line $y=x$ becomes as flat as the horizon --- and then using the persistence diagram points to trace out the peaks in a mountain landscape profile.
A landscape can then be discretized by taking a finite sample of the function values, allowing it to be used in machine learning tasks: see for instance~\cite{kovacev-nikolic2016protein_binding}.
From barcodes,~\citet{PersistenceImages} create \emph{persistence images}, a Euclidean vectorization enabling a diverse class of machine learning tools to be applied (see also~\citet{chen2015statistical,7299106}).
A persistence image is created by taking a sum of Gaussians, one centered on each point in a persistence diagram, and then pixelating that surface to form an image.
By analogy, recall that in point cloud persistent homology, one ``blurs their vision" when looking at a dataset by replacing each data point with a ball --- this is similar to the process of ``blurring one's vision" when looking at a persistence diagram in order to create a persistence image.

Persistence landscapes were defined as part of an effort to give a firm statistical foundation to persistent homology.
In fact, \citet{bubenik2015statistical} proves a strong law of large numbers and a central limit theorem for persistence landscapes.
This allows one to discuss hypothesis testing with persistent homology.
Another approach to hypothesis testing is given by~\citet{robinson2017hypothesis}.
Other statistical approaches include~\citet{Fasy2014}, which describes confidence intervals and a statistical approach to distinguishing important features from noise,~\citet{divol2019choice,JMLR:v20:18-618}, which consider probability density functions for persistence diagrams, and \citet{Maroulas2020Bayesian}, which describes a Bayesian framework.
See~\citet{wasserman2018topological} for a review of statistical techniques in the context of topological data analysis.

Persistence landscapes and images are only two of the many different methods that have recently been invented in order to transform persistence barcodes into machine learning input.
Algorithms that require only a distance matrix, such as many clustering or dimensionality reduction algorithms, can be applied on the bottleneck or Wasserstein distances between persistence barcodes~\citep{cohen2007stability,mileyko2011probability,kerber2017geometry}.
Other techniques for vectorizing persistence barcodes involve heat kernels~\citep{carriere2015stable}, entropy~\citep{merelli2015topological,atienza2018stability}, rings of algebraic functions~\citep{adcock2016ring}, tropical coordinates~\citep{kalivsnik2019tropical}, complex polynomials~\citep{di2015comparing}, and optimal transport~\citep{carriere2017sliced}, among others.
Some of these techniques, including those by~\citet{zhao2019learning} and~\cite{divol2019choice}, allow one to learn the vectorization parameters that are best suited for a machine learning task on a given dataset.
Others allow one to plug persistent homology information directly into a neural network~\citep{hofer2017deep}.
Recent research on incorporating persistence as input for machine learning is vast and varied, and the above collection of references is far from complete.

As for the creation of new algorithms, persistent homology has recently been applied to regularization, a technique used in machine learning that penalizes overly complicated models to avoid overfitting. 
\citet{pmlr-v89-chen19g} propose a ``topological penalty function" for classification algorithms, which encourages a topologically simple decision boundary.  
Their method is based on measuring the relative importance of various connected components of the decision boundary via 0-dimensional persistent homology.
They show how the gradient of such a penalty function can be computed, which is important for use in machine learning algorithms, and demonstrate their method on several examples.
Similar work using topological methods to examine a decision boundary can also be found in~\citet{varshney2015decisionboundaries} and~\citet{pmlr-v97-ramamurthy19a}.

Finally, other recent work has used persistent homology to analyze neural networks.
~\citet{JMLR:v21:20-345} provide experimental evidence that neural networks operate by simplifying the topology of a dataset.
They examine the topology of a dataset and its images at the various layers of a neural network performing classification, finding that the corresponding barcodes become simpler as the data progresses though the network.  
Additionally, they observe the effects of different shapes of neural networks and different activation functions.  
They find that deeper neural networks have a tendency to simplify the topology of a dataset more gradually than shallow networks, and that networks with ReLU activation tend to simplify topology more in the earlier layers of a network than other activation functions.

\section{Conclusion}

Topological tools are often described as being able to stitch local data together in order to describe global features: from local to global.
The history of applied topology, however, has in some sense gone in the reverse direction --- from global to local --- as surveyed above!
Applied topology was developed in part to summarize global features in a point cloud dataset, as in the examples of the conformations of the cyclo-octane molecule or the collection of $3\times 3$ pixel patches from images.
If global shapes are the focus, long persistent homology bars are interpreted as the relevant features, while small bars are often disregarded as sampling artifacts or noise.
However, in more recent applications, and in particular when using applied topology in concert with machine learning, it is often many short persistent homology bars that \emph{together} form the signal.
One of the biggest benefits of applied topology is that one need not choose a scale beforehand: persistent homology provides a useful summary of both the local and global features in a dataset, and this summary has been made accessible for use in machine learning tasks.

We have seen how the short bars can be a measure of local geometry, texture, curvature, and fractal dimension; their sensitivity to various features of datasets leads to the wide variety of applications surveyed here.
Because persistent homology provides a concise, reductive view of the geometry of a dataset, for instance in the examples studying brain artery trees or hexagonal grids, it is not hard to imagine the potential applications to machine learning problems.  
This has led to recent techniques that turn barcodes into machine learning input, exemplified by persistence landscapes and persistence images. 
We hope that this wealth of recent work, which has shifted more attention to short persistent homology bars and the geometric information they summarize, will inspire further research at the intersection of applied topology, local geometry, and machine learning.





\section*{Funding}
This material is based upon work supported by the National Science Foundation under Grant Number 1934725.

\bibliographystyle{plainnat}
\bibliography{LocalToGlobal}

\end{document}